\documentclass[12pt]{amsart}
\usepackage{amsfonts}
\usepackage{amsthm}
\usepackage{enumerate}
\usepackage{helvet,courier}
\usepackage{mathptmx,amsmath}
\usepackage{type1cm}
\DeclareMathAlphabet{\mathcal}{OMS}{cmsy}{m}{n}
\DeclareMathAlphabet{\mathbbold}{U}{bbold}{m}{n}  
\usepackage{amsthm}
\usepackage{graphicx}       
\usepackage{multicol}
\usepackage{mathrsfs}
\usepackage{amssymb}
\usepackage{verbatim}
\usepackage{esint}

\usepackage{mathtools}
\usepackage{tikz}
\usepackage[hidelinks]{hyperref}

\usepackage{pgfplots}
\usepgfplotslibrary{fillbetween}
\pgfplotsset{compat=1.18}

\usetikzlibrary{arrows.meta}

\usepackage[dvipsnames]{xcolor}

\theoremstyle{plain}
\newtheorem{thm}{Theorem}
\newtheorem{lm}{Lemma}[section]
\newtheorem{cor}[lm]{Corollary}
\newtheorem{prop}[lm]{Proposition}

\newtheorem{conj}{Conjecture}

\theoremstyle{remark}
\newtheorem{rmk}{Remark}

\theoremstyle{definition}
\newtheorem{defn}{Definition}

\newcommand{\bnu}{\begin{enumerate}}
\newcommand{\enu}{\end{enumerate}}
\newcommand{\bpf}{\begin{proof}}
\newcommand{\epf}{\end{proof}}

\newcommand{\sset}{\subset}

\newcommand{\al}{\alpha}
\newcommand{\be}{\beta}

\newcommand{\La}{\Lambda}

\newcommand{\ep}{\epsilon}

\newcommand{\bbz}{\mathbb{Z}}

\newcommand{\bbr}{\mathbb{R}}

\newcommand{\bbrn}{\mathbb R^n}

\newcommand{\bbc}{\mathbb{C}}

\newcommand{\ges}{\gtrsim}
\newcommand{\les}{\lesssim}

\newcommand{\f}{\frac}

\newcommand{\tf}{\tfrac}

\begin{document}

\author{Sian Fang}
\address{Sian Fang, School of Mathematical Sciences,
Fudan University, People's Republic of China}
\email{24110180013@m.fudan.edu.cn}

\author{Danqing He}
\address{Danqing He, School of Mathematical Sciences,
Fudan University, People's Republic of China}
\email{hedanqing@fudan.edu.cn}

\author{Xinyu Zhu}
\address{Xinyu Zhu, School of Mathematical Sciences,
Fudan University, People's Republic of China}
\email{23110180063@m.fudan.edu.cn}

\thanks{All authors are supported by  National Key R$\&$D Program of China (No. 2021YFA1002500), and NNSF of China (No. 12322105). In addition, D. He is supported by the New Cornerstone Science Foundation. }

\title{Products of Translates: Bounds and Constructions}
\date{}
\maketitle

\begin{abstract}
    We study products of translates of the form $(A+d)(A+d)$ with $d\in D$, where both $A$ and $D$ are finite subsets of $\bbr$. We show that
    $|(A+a)(A+a)|\ge C|A|^{5/3}$ for some $a\in A$, and that $|(A+D)(A+D)|\le N^{2-c(\al)}$ for some $c(\al)>0$ when $|D|=|A|^\alpha$ with $\al<1$.
\end{abstract}


\section{Introduction}

Let $A,B\subset\bbr$ be finite nonempty sets.
Their sumset is
$
A+B:=\{a+b:a\in A,\, b\in B\}.
$
If $B=\{x\}$, we write $A+B=A+x$ and call this set a translate of $A$. Similarly, define
$
AB:=\{ab:a\in A,\, b\in B\},
$
and, if $0\notin B$,
$
A/B:=\{a/b:a\in A,\, b\in B\}.
$
The sum-product conjecture raised by Erd\H{o}s--Szemer\'edi \cite{E83} asks whether every finite nonempty set $A\subset\bbz$ satisfies
\[
\max\{|A+A|,|AA|\}\gtrsim_\epsilon |A|^{2-\epsilon}
\qquad\text{for every }\epsilon>0.
\]
The analogous question for subsets of $\bbr$ was also widely studied. The sum-product conjecture expresses the expectation that strong additive and multiplicative structure cannot coexist. Arithmetic progressions and positive geometric progressions illustrate the contrasting  behaviors  of the two structures.
Elekes \cite{E97} proved that 
$$
\max\{|A+A|,|AA|\}\gtrsim |A|^{5/4}
$$
via the Szemer\'edi-Trotter theorem \cite{S83}.
An elegant argument introduced by  Solymosi  \cite{S09}
yields the exponent $4/3$, up to a logarithmic factor. 
More recently, Bloom \cite{B25} improve the exponent to $\frac{1270}{951}$.

The distinction between the integers and the reals is now essential. Bloom, Sawin, Schildkraut and Zhelezov \cite{B26} constructed arbitrarily large real sets for which
\[
\max\{|A+A|,|AA|\}\leq |A|^{2-c}
\]
with an absolute $c>0$, disproving the real-number version of the sum-product conjecture.

The sum-product phenomenon also leads to questions involving mixed operations. One example is 
\[
(A+a)A=\{(x+a)y:x,y\in A\}.
\]
For $a=0$, this is the product set $AA$, which may be small. For $a\neq0$, one expects substantially greater expansion. 
\cite{G10} proved the  estimate
$
|A(A+1)|\gtrsim |A|^{5/4},
$
where 
\[
\max\{|AA|,|(A+1)(A+1)|\}\gtrsim |A|^{5/4}
\]
is also proved implicitly.
These results were later improved by 
\cite{J13}, \cite{H19}, and \cite{S22}.
These results quantify the expansion produced by multiplying a set by a nontrivial translate of itself.

Recent work also gives counterexamples to near-quadratic growth for fixed translates. 
Adapting the construction in \cite{B26}, 
\cite{R26} constructed arbitrarily large real sets such that, for an absolute $c>0$,
$
\max\{|AA|,|A(A+1)|,|(A+1)(A+1)|\}
\lesssim N^{2-c}.
$
\cite{R26} showed that for any fixed $k$, there exists a constant $c(k)$ such that
\begin{equation}\label{e050}
    \max_{d\in\{0,1,\dots,k\}}|(A+d)(A+d)|\lesssim N^{2-c(k)} .
\end{equation}

In \eqref{e050}, $k$ is a fixed number, so the behavior of $   \max_{d\in\{0,1,\dots,k\}}|(A+d)(A+d)|$ is expected to be similar to the case $k=1$ as $N\to\infty$. But what do we have if $k$ depends on $N$? 
A natural model motivated by this question is  
$$\max_{a\in A}|(A+a)(A+a)|.$$
This model can be regarded as the diagonal version of 
$\max_{a,b\in A}|(A+a)(A+b)|$
studied by 
Roche-Newton \cite{R15},
who proved that there exist $a,b\in A$ such that
\[
|(A+a)(A+b)|\gtrsim\frac{|A|^2}{\log(2|A|)},
\]
by combining energy estimates with the Szemer\'edi--Trotter theorem.

To formulate our problems, let us give some definitions first.
Let  \(A\subset\mathbb R\) be a finite  nonempty set with  \(N=|A|\). We
use the notation
\[
 F(d)=(A+d)(A+d)=\{(x+d)(y+d):x,y\in A\},\qquad f(d)=|F(d)|.
\]
Our main questions concern (i) the size of $F(d)$ when $d$ varies over $\bbr$, (ii) the number of translations for which $F(d)$ is small, and, especially, (iii) the existence of a large product set when the translation is restricted to a finite set $D$ with $|D|\sim N^{\alpha}$ for some $\alpha\ge0$.

Commutativity gives the universal upper bound
\[
f(d)\leq \frac{N(N+1)}2.
\]
The bound fails to be attained for at most $N^4$ many $d$, as shown by Proposition \ref{t40}.
This result gives partial answers to questions (i) and (ii), whose  proof is given in Section~\ref{s40}.
A concrete illustration is the geometric progression $A=\{2^k:1\leq k\leq N\}$, for which
\begin{equation}\label{e15}
|(A+d)(A+d)|=\frac{N(N+1)}2
\end{equation}
whenever $0<|d|<1$; see Section~\ref{exam6.1}. The examples in Section~\ref{exam6.1} further show that this geometric progression has $|(A+a)(A+a)|\gtrsim N^2$ for every $a\neq 0$, despite having a product set of only linear size at $a=0$.

\medskip

However, the above generic observation does not settle the pinned problem: all $N$ admissible shifts may lie among the exceptional translations without any additional information. Concerning the question (iii), we raise the following conjecture.

\begin{conj}\label{t31}

For every $\epsilon>0$, there exists a positive constant $C_\epsilon$ such that, for any $D\subset\bbr$ with $|D|\sim N$, 
$$
\sup_{d\in D}|(A+d) (A+d)|\ge C_\ep |A|^{2-\ep}
$$
holds for any finite nonempty set $A\sset\bbr$.

\end{conj}

\begin{rmk}
    One may consider also a weak version
    \begin{equation}\label{e057}
        M(A)\gtrsim_\epsilon |A|^{2-\epsilon}.
    \end{equation} But the methods used in this note apply to all  $D$ which do not have to be $A$.
\end{rmk}

To sutdy Conjecture~\ref{t31}, we have the following definitions.
For a finite nonempty  set \(D\subset\mathbb R\), put
\begin{equation}\label{e91}
 M_D(A)=\max_{d\in D}f(d),\qquad M(A)=M_A(A).
\end{equation}

\begin{defn}\label{d901}
For \(\alpha\ge0\), let \(\beta(\alpha)\) be the supremum of all
\(b\in[0,2]\) with the following property. For every
\(\varepsilon>0\) and every \(\kappa\ge1\), there exist
\(c=c(b,\alpha,\varepsilon,\kappa)>0\) and
\(N_0=N_0(b,\alpha,\varepsilon,\kappa)\) such that
\begin{equation}\label{e90}
 M_D(A)\ge cN^{b-\varepsilon}
 \quad\text{if}\quad
 |A|=N\ge N_0,\qquad
 \kappa^{-1}N^\alpha\le |D|\le\kappa N^\alpha.
\end{equation}
The constants are uniform over all such real sets \(A,D\).
\end{defn}

Roughly speaking, $\be(\al)$ is the largest number such that
$
M_D(A)\gtrapprox N^{\be(\al)}
$
for all $A$ with $|A|=N$ and  $D$ with $|D|\sim N^{\al}$. Conjecture~\ref{t31} claims that  $\beta(1)=2$.

Using the main ingredients in \cite{B26} and \cite{R26} , we obtain the following counterexample, which paritially supported the rationality of Conjecture~\ref{t31}.

\begin{thm}
\label{t30}
There is an absolute constant \(c_*>0\) such that for every fixed
\(0<\alpha<1\), there exist arbitrarily large finite sets \(A\), with $N:=|A|$,  and
corresponding finite sets \(D\) such that
\[
 |D|\sim N^\alpha,\qquad
 M_D(A)\le N^{2-c_*(1-\alpha)}.
\]
Moreover, our construction gives
$$
|(A+D)(A+D)|\le N^{2-c_*(1-\al)}.
$$
\end{thm}

The multiplicative energy of $A$ and $B$ is defined as
$$
E_{\times}(A,B)
:=\#\{(a_1,a_2,b_1,b_2)\in A^2\times B^2:
             a_1b_1=a_2b_2\}.
$$
Concerning \eqref{e057}, the weak version of Conjecture~\ref{t31}, we have the following result.

\begin{thm}
\label{t32}
Let $A\subset\bbr$ be finite and nonempty, with $|A|=N$.
\begin{enumerate}
\renewcommand{\labelenumi}{(\roman{enumi})}
\item The averaged multiplicative energy satisfies
$$
\sum_{a\in A}E_\times(A+a)\lesssim |A|^{10/3}.
$$
\item There exists $a\in A$ such that $|(A+a)(A+a)|\gtrsim |A|^{5/3}$.
\end{enumerate}
\end{thm}

We can actually prove a   result stronger than Theorem~\ref{t32}.
Let
\begin{equation}\label{e94}
 L(\alpha)=
 \begin{cases}
 1,&\alpha=0,\\
 \tfrac{49}{38},&0<\alpha\le\tfrac{11}{38},\\
 1+\alpha,&\tfrac{11}{38}\le\alpha\le\tfrac12,\\
 \tfrac43+\tfrac\alpha3,&\tfrac12\le\alpha\le2,\\
 2,&\alpha\ge2.
 \end{cases}
\end{equation}

\begin{thm}\label{t58}
    For all \(\alpha\ge0\), we have
\begin{equation}\label{e95}
 L(\alpha)\le\beta(\alpha).
\end{equation}
\end{thm}

\begin{figure}[!htbp]
\centering
\usetikzlibrary{patterns}
\resizebox{0.97\linewidth}{!}{\begingroup
\definecolor{bcBlue}{HTML}{17639B}
\definecolor{bcOrange}{HTML}{C56339}
\definecolor{bcPurple}{HTML}{76518D}
\definecolor{bcDark}{HTML}{263545}
\definecolor{bcMuted}{HTML}{64748B}
\definecolor{bcBlueFill}{HTML}{DDECF9}
\definecolor{bcGreenFill}{HTML}{E1F1DE}
\definecolor{bcOrangeFill}{HTML}{FAE4D6}
\definecolor{bcGrid}{HTML}{E6EAEF}
\begin{tikzpicture}[x=4.16cm,y=4.55cm,
    font=\sffamily\small]
  \path[use as bounding box] (-0.41,0.57) rectangle (3.16,2.115);

  \begin{scope}
    \clip (0,0.99) rectangle (3.1,2.075);
    \path[fill=bcBlueFill]
      (0,1) -- (0,1.289473684)
      -- (0.289473684,1.289473684)
      -- (0.5,1.5) -- (2,2) -- (3.1,2) -- (3.1,1) -- cycle;
    \path[fill=bcGreenFill]
      (0,1.289473684) -- (0.289473684,1.289473684)
      -- (0.5,1.5) -- (2,2) -- (1,2) -- (0,1.72) -- cycle;
    \path[fill=bcOrangeFill]
      (0,1.72) -- (1,2) -- (0,2) -- cycle;
    \path[pattern=north east lines,pattern color=bcOrange,opacity=.24]
      (0,1.72) -- (1,2) -- (0,2) -- cycle;
    \foreach \yy in {1,1.289473684,1.5,1.666666667,2}{
      \draw[bcGrid,thin] (0,\yy) -- (3.1,\yy);
    }
    \draw[bcMuted,densely dotted,thin] (0,2) -- (1,2);
    \draw[bcMuted,densely dotted,thin] (1,0.99) -- (1,2.075);
    \draw[bcMuted,densely dotted,thin] (2,0.99) -- (2,2.075);
    \draw[bcBlue,line width=2.1pt]
      (0,1.289473684) -- (0.289473684,1.289473684)
      -- (0.5,1.5) -- (2,2) -- (3.1,2);
    \draw[bcOrange,line width=2pt,dash pattern=on 8pt off 5pt]
      (0,1.72) -- (1,2);
    \draw[bcPurple,line width=2.3pt,dash pattern=on 2pt off 5pt]
      (1,2) -- (2,2);
    \draw[bcBlue,line width=2.45pt] (2,2) -- (3.1,2);
  \end{scope}

  \draw[bcDark,line width=.75pt] (-.025,.99) -- (3.1,.99);
  \draw[bcDark,line width=.75pt] (-.025,.99) -- (-.025,2.075);
  \foreach \xx in {0,.5,1,2,3}{
    \draw[bcDark,line width=.65pt] (\xx,.99) -- (\xx,.975);
  }
  \foreach \yy in {1,1.289473684,1.5,1.666666667,2}{
    \draw[bcDark,line width=.65pt] (-.025,\yy) -- (-.043,\yy);
  }

  \node[anchor=north,text=bcDark] at (0,.96) {$0$};
  \node[anchor=north,text=bcDark] at (.5,.96) {$\frac{1}{2}$};
  \node[anchor=north,text=bcDark] at (1,.96) {$1$};
  \node[anchor=north,text=bcDark] at (2,.96) {$2$};
  \node[anchor=east,text=bcDark] at (-.055,1) {$1$};
  \node[anchor=east,text=bcDark] at (-.055,1.289473684) {$49/38$};
  \node[anchor=east,text=bcDark] at (-.055,1.5) {$3/2$};
  \node[anchor=east,text=bcDark] at (-.055,1.666666667) {$5/3$};
  \node[anchor=east,text=bcDark] at (-.055,2) {$2$};
  \node[text=bcDark,font=\sffamily\large] at (1.52,.84) {$\alpha$};
  \node[text=bcDark,font=\sffamily\large,rotate=90] at (-.34,1.57) {$\beta(\alpha)$};

  \node[text=bcDark,font=\sffamily\footnotesize] at (.34,1.93) {Counterexamples};
  \node[text=bcDark] at (1.43,1.93) {Unknown};
  \node[text=bcDark] at (.67,1.18) {Positive results};

  \fill[bcBlue] (0,1) circle[radius=1.8pt];
  \fill[white] (0,1.289473684) circle[radius=2.3pt];
  \draw[bcBlue,line width=1.35pt] (0,1.289473684) circle[radius=2.3pt];
  \fill[white] (0,1.72) circle[radius=2.3pt];
  \draw[bcOrange,line width=1.35pt] (0,1.72) circle[radius=2.3pt];

  \node[anchor=west,text=bcDark,font=\small] at (2.57,1.055) {$\alpha\to\infty$};
  \draw[bcDark,-stealth,line width=.75pt] (2.90,1.04) -- (3.08,1.04);

  \path[fill=bcBlueFill,draw=bcBlue,line width=.6pt]
    (.61,.72) rectangle (.74,.747);
  \node[anchor=west,text=black,font=\sffamily\scriptsize]
    at (.775,.734) {Proved};
  \path[fill=bcGreenFill,draw=green!55!black,line width=.6pt]
    (.61,.655) rectangle (.74,.682);
  \node[anchor=west,text=black,font=\sffamily\scriptsize]
    at (.775,.668) {Unknown};
  \path[fill=bcOrangeFill,draw=bcOrange,line width=.6pt]
    (1.12,.72) rectangle (1.25,.747);
  \path[pattern=north east lines,pattern color=bcOrange]
    (1.12,.72) rectangle (1.25,.747);
  \node[anchor=west,text=black,font=\sffamily\scriptsize]
    at (1.285,.734) {Excluded};
  \draw[bcBlue,line width=2pt] (1.12,.668) -- (1.25,.668);
  \node[anchor=west,text=black,font=\sffamily\scriptsize]
    at (1.285,.668) {$L(\alpha)$};
  \draw[bcOrange,line width=1.8pt,dash pattern=on 7pt off 4pt]
    (1.73,.734) -- (1.86,.734);
  \node[anchor=west,text=black,font=\sffamily\scriptsize]
    at (1.90,.734) {Counterexample bound};
  \draw[bcPurple,line width=1.9pt,dash pattern=on 2pt off 4pt]
    (1.73,.668) -- (1.86,.668);
  \node[anchor=west,text=black,font=\sffamily\scriptsize]
    at (1.90,.668) {Conjectured plateau};

\end{tikzpicture}
\endgroup
}
\caption{Known lower and upper bounds for $\beta(\alpha)$. The visible slope for  counterexamples does not specify the proved constant $c_*$}
\label{f0001}
\end{figure}
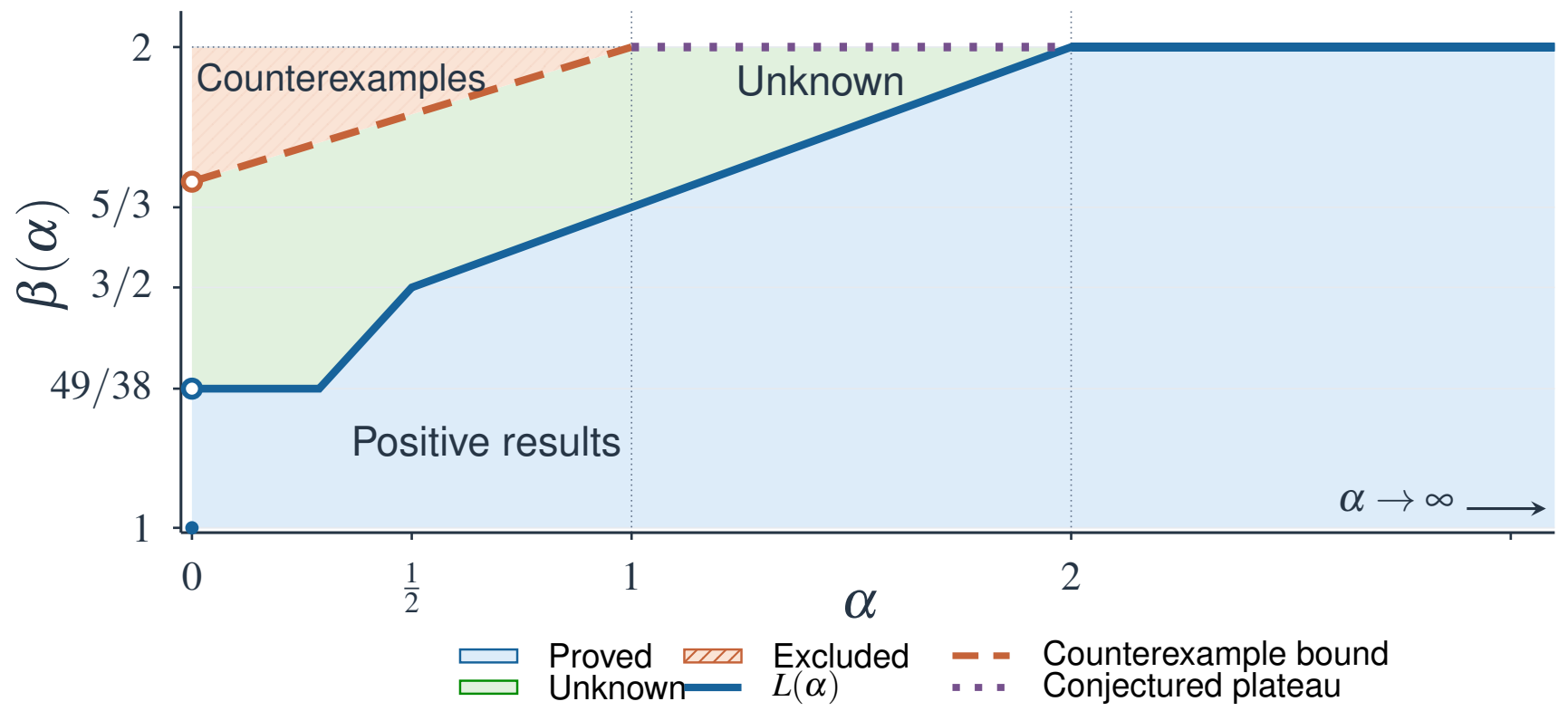

We explain the main idea of the proof of Theorem~\ref{t58} by dicussing the case $\al=1$ and $D=-A$ briefly.
We first express the energy equations for $A-a$ as collinearity relations between grid points and a diagonal pin $(a,a)$. 
We reduce the problem to counting the number of  $\mathcal L_a^r$ with $a\in A$, which consist of 
$r$-rich lines  (with respect to $A\times A$) through $(a,a)$.
Theorem~\ref{t51} yields
\[
|\bigcup_{a\in A}\mathcal L_a^r|\lesssim
\min\left\{\frac{N^3}{r},\frac{N^{7/2}}{r^{5/2}}\right\},
\]
whose proof uses the ordering of $A$: many consecutive points on a rich line have small gaps in their coordinate indices. Reorganizing these pairs as collinear triples reduces the count to a second-moment incidence estimate on a smaller subset of the grid $A\times A$, which is bounded by using the Szemer\'edi-Trotter theorem and a bush argument. From the analytic viewpoint, this brings us some orthogonality, which leads to the final improvement. The two terms in the improved bound balance at $r\sim N^{1/3}$, and summation over dyadic richness levels produces the exponent $10/3$ in the averaged energy.

This paper is organized as follows.
In Section~\ref{s092} we prove the main theorem.
In Section~\ref{bush} we prove the main technical result using a bush argument.
Section~\ref{s40} gives exceptional-set estimates.
Section~\ref{s9} presents some counterexamples.

\medskip

\noindent{\bf Notations.} Throughout, $U\lesssim V$ means $U\leq CV$ for an absolute constant $C$, and $U\sim V$ means that both $U\lesssim V$ and $V\lesssim U$ hold. A subscript indicates allowed dependence of the constant. Also, we write $U\lessapprox V$ if $U\lesssim_{\ep} N^{\ep}V$ for any $\ep>0$.  
Unindexed sums of representation functions run over their finite supports. 

\medskip

\noindent{\bf Use of AI.} The authors formulated the main conjecture and developed the proof strategy. AI contributed an essential idea to the lower-bound proof and proposed a construction yielding small product sets. The authors independently checked every proof detail and rewrote the manuscript, taking full responsibility.

\section{Lower bounds via energy}
\label{s092}

We prove Theorem~\ref{t58} in this section.
We utilize a classical approach to reduce the problem to an incidence estimate via a intermediate multiplicative energy control. See, for instance, \cite{R15} for an argument of this type.

For finite sets $A,B\subset\bbr$, define
$
r_{AB}(z):=\#\{(x,y)\in A\times B:xy=z\}
$
and
\[
\begin{aligned}
E_{\times}(A,B)
&:=\#\{(a_1,a_2,b_1,b_2)\in A^2\times B^2:
             a_1b_1=a_2b_2\}=\sum_{z\in AB}r_{AB}(z)^2.
\end{aligned}
\]
We write $E_{\times}(A)=E_{\times}(A,A)$. 
The Cauchy--Schwarz inequality gives
\begin{equation}\label{CS}
E_{\times}(A,B)|AB|\geq |A|^2|B|^2.
\end{equation}

\begin{lm}
    Let $D\subset \mathbb R$ be a finite nonempty set. Then
    \begin{equation}\label{e054}
        M_D(A)\sum_{d\in D}E_\times (A+d)\ge N^4|D|.
    \end{equation}
\end{lm}

\bpf

It follows from \eqref{CS} that
$$
\sum_{d\in D} [f(d)]^{-1}\le N^{-4}\sum_{d\in D} E_\times(A+d),
$$
which yields \eqref{e054}.

\epf

For each $d\in D$, the number of ordered pairs $(x,y)\in A^2$ with $(x+d)(y+d)=0$ is $2N-1$. Thus the total zero-product contribution to $\sum_{d\in D}E_\times(A+d)$ is bounded by $O(N^3)$. The nonzero contribution is
\begin{equation}
\begin{aligned}
\sum_{d\in D}\#\left\{(x,y,z,w)\in A^4:
\frac{y+d}{z+d}=\frac{w+d}{x+d},\quad x,y,z,w\neq -d\right\}.
\end{aligned}
\label{e81}
\end{equation}
The quotient equality says that $(z,y)$ and $(x,w)$ lie on the same line through $(-d,-d)$.

Fix finite subsets $A $ and $X$ of $\bbr$.
Write $\Delta=\{(u,u):u\in\bbr\}$ and define, for $x\in X$,
\[
\mathcal L_{A,x}:=\{\ell\neq\Delta:(x,x)\in\ell,\ 
\ell\cap(A\setminus\{x\})^2\neq\emptyset\}.
\]
Here and below $\ell$ denotes an affine line. These lines are neither horizontal nor vertical. 
Every line in $\mathcal L_{A,d}$ contains at least two grid points in $A\times A$.
The diagonal is excluded because it belongs to every bush, which is a family of lines through one fixed point. 
For an integer $r$ with $2\leq r\leq N$, put
\[
\mathcal L_{A,x}^r:=\{\ell\in\mathcal L_{A,x}:r\leq|\ell\cap(A\times A)|<2r\},
\qquad L_{A,X}^r:=\bigcup_{x\in X}\mathcal L_{A,x}^r,
\qquad L_{A,X}:=\bigcup_{x\in X}\mathcal L_{A,x}.
\]

We reduce the control of energy to an incidence estimate, as revealed by the next result.

\begin{lm}\label{IC}

\begin{equation}\label{e82}
\sum_{d\in D}E_\times(A+d)
\lesssim N^3+N^2|D|+
\sum_{\text{dyadic }r\in[16,N]}r^2|L_{A,-D}^r|.
\end{equation}

\end{lm}
\bpf

Separating the diagonal contribution gives 
\begin{equation}
\begin{aligned}
&\#\left\{(x,y,z,w)\in A^4:
\frac{y+d}{z+d}=\frac{w+d}{x+d},\quad x,y,z,w\neq -d\right\}\\
&\le N^2+\sum_{\ell\in\mathcal L_{A,-d}}
|\ell\cap(A\times A)|^2.
\end{aligned}
\label{e83}
\end{equation}

We observe that distinct bushes are disjoint since a line containing two distinct diagonal pins would be $\Delta$. 

In particular, $\sum_{d\in D}|\mathcal L_{A,-d}^r|=|L_{A,-D}^r|$. 
Within each bush $\mathcal L_{A,-d}$, $\ell\cap A\times A$ are disjoint except for $(-d,-d)$. This implies that 
$$
\sum_{\ell\in\mathcal L_{A,-d}^r: r\le 16}
|\ell\cap(A\times A)|^2\les N^2.
$$

Summing \eqref{e83} over dyadic richness levels gives
\[
\sum_{d\in D}E_\times(A+d)
\lesssim N^3+N^2|D|+\sum_{\text{dyadic } r\in [16,N]}r^2|L_{A,-D}^r|.
\]

\epf

We  record  two elementary bounds.

\begin{lm}
\label{t50}
For every integer $r$ with $2\leq r\leq N$,
\begin{equation}
\left|L_{A,-D}^r\right|
\leq
C\min\left\{\frac{N^2|D|}{r},\frac{N^4}{r^3}\right\}.
\label{e84}
\end{equation}
\end{lm}

\begin{proof}

For a fixed $d\in D$, the sets $\ell\cap(A\times A\setminus\Delta)$, with $\ell\in\mathcal L_{A,-d}^r$, are disjoint. Consequently,
\[
(r-1)|\mathcal L_{A,-d}^r|\leq(N-1)^2,
\qquad
\sum_{d\in D}|\mathcal L_{A,-d}^r|\lesssim\frac{N^2|D|}{r}.
\]

The Szemer\'edi--Trotter theorem about $r$-rich lines says that
$$
\left|L_{A,D}^r\right|
\lesssim
\frac{|A\times A|^2}{r^3}
+
\frac{|A\times A|}{r}
\lesssim
\frac{N^4}{r^3}
$$
since $r\leq N$.
\end{proof}

We can actually obtain the following improvement over \eqref{e84}.

\begin{thm}\label{BA}
\label{t51}
For every integer $r$ with $2\leq r\leq N$, and every real number $\alpha$ with $\alpha>0$, if  a subset $D$ of $\bbr$ satisfies $|D|\sim N^\alpha$, we can get that
$$
\left|L_{A,D}^r\right|
\lesssim
\min\left\{\frac{N^{2+\alpha}}{r},\frac{N^{3+\frac\al 2}}{r^{5/2}}+\frac{N^3}{r^2}\right\}.
$$
\end{thm}

We  invoke the following result of Stevens and Warren \cite{S22} to attack $\be(\al)$ for small $\al$. 
\begin{lm}\label{t52}
    There is an absolute constant $C>0$ such that for all finite $A\subset\bbr$ with $|A|=N\ge2$,
    $$
    |AA|+|(A+1)(A+1)|\ges N^{49/38}(\log N)^{-C}.
    $$
\end{lm}
The above estimates have an immediate corollary for the product of two distinct shifts.
\begin{cor}\label{t53}
     There is an absolute constant $C>0$ such that for all finite $A\subset\bbr$ with $|A|=N\ge2$ and $d_0\ne d_1\in\bbr$,
    \begin{equation}\label{e88}
    |(A+d_0)(A+d_0)|+|(A+d_1)(A+d_1)|\ges N^{49/38}(\log N)^{-C}.
    \end{equation}
    Here the implicit constant is independent of $d_0$ and $d_1$. In particular,
    $$
    M_D(A)\gtrapprox N^{49/38}
    $$
    for all $D\subset\bbr$ with $|D|\ge2$.
\end{cor}
\begin{proof}
    As $d_0\ne d_1$, we can define $B=\frac{A+d_0}{d_1-d_0}$. Notice that
    $$(A+d_1)(A+d_1)=(d_1-d_0)^2(B+1)(B+1).$$
    Then the left hand side of \eqref{e88} is bounded by $|BB|+|(B+1)(B+1)|$. Then Lemma~\ref{t52} gives the desired bounds.
\end{proof}

We are ready to prove Theorem~\ref{t58}. Let us first recall
$$
 L(\alpha)=
 \begin{cases}
 1,&\alpha=0,\\
 \tfrac{49}{38},&0<\alpha\le\tfrac{11}{38},\\
 1+\alpha,&\tfrac{11}{38}\le\alpha\le\tfrac12,\\
 \tfrac43+\tfrac\alpha3,&\tfrac12\le\alpha\le2,\\
 2,&\alpha\ge2.
 \end{cases}
$$
Our goal is to show that $\be(\al)\ge L(\al)$.
\begin{proof}[Proof of Theorem~\ref{t58}]
The case $\al=0$ is trivial. We may assume $\al>0$ below.

It follows from Theorem~\ref{BA} that
\begin{align*}
    \sum_{\substack{r=2^j\\16\leq r\leq N}}r^2|L_{A,-D}^r|&\les\sum_{\substack{r=2^j\\16\leq r\leq N^{2/3-\al/3}}}rN^{2+\al}+\sum_{\substack{r=2^j\\\max\{16,N^{2/3-\al/3}\}\leq r\leq N}}(\f{N^{3+\f{\al}2}}{r^{1/2}}+N^3)\\
    &\les N^{8/3+2\al/3}+N^3 \log N.
\end{align*}
Applying Lemma~\ref{IC}, we obtain
\begin{equation}\label{e056}
\sum_{d\in D}E_\times(A+d)
\lesssim N^{8/3+2\al/3}+N^{2+\al}+N^3\log(2N),
\end{equation}
    which combined with \eqref{e054} yields that 
    \begin{equation}\label{MDA}
    N^{\al+4}\les M_D(A)(N^{\f{2}{3}\al+\f{8}{3}}+N^{\al+2}+N^3\log N).
    \end{equation}
    The case $0<\al\le 11/38$ is a consequence of Corollary~\ref{t53}. All other cases follow from \eqref{MDA} by a routine calculation.

    In summary, we have $\be(\al)\ge L(\al)$ as desired.
\end{proof}

\begin{proof}[Proof of Theorem~\ref{t32}]
(i) follows from \eqref{e056}. (ii) follows from Theorem~\ref{t58} when $\al=1$.

\end{proof}

\begin{rmk}
The same conclusions hold for $(A+a)(A+ka)$, where $k\in\bbr\setminus\{0\}$, with absolute constants independent of $k$. Indeed, set $B=A\cup k^{-1}A$. For every $a\in A\subset B$,
\[
E_\times(A+a,A+ka)
=E_\times(A+a,k^{-1}A+a)
\leq E_\times(B+a).
\]
Summing and applying Theorem~\ref{t32} to $B$, with $|B|\leq2N$, gives the energy bound $\lesssim N^{10/3}$; inequality~\eqref{CS} gives the product-set bound. This argument requires $k\neq0$ and gives no corresponding assertion for the degenerate case $k=0$, which was studied in \cite{M15} and \cite{B25}.
\end{rmk}

\section{Improved Rich-Line bounds via bushes}\label{bush}

The proof of Theorem~\ref{t51} uses index proximity, as in related arguments of \cite{B22,S24,C26}.

\begin{proof}[Proof of Theorem~\ref{t51}]
For $\al\ge 2$, because
$$
\frac{N^4}{r^3}\lesssim \min\left\{\frac{N^{2+\alpha}}{r},\frac{N^{3+\frac\al 2}}{r^{5/2}}+\frac{N^3}{r^2}\right\},
$$
we can finish the proof from Lemma~\ref{t50}. So we only need ton consider the case $0<\al<2$.

It suffices to prove the second bound. For $2\leq r<16$, it follows from the first bound after changing the absolute constant. We may therefore assume $r\geq16$.

Let us order $A$ by
$$
a_1<a_2<\cdots<a_N.
$$
We define $I(a_i)=i$, the index of $a_i$ in $A$, and the index distance
$$
d_I(a_i,a_j)=\left|I(a_i)-I(a_j)\right|.
$$
For an integer $N>K\geq1$, define
$$
P_K:=\left\{(b,c)\in A\times A:0<d_I(b,c)\leq K\right\}.
$$
Notice that
\begin{equation}
|P_K|=2\sum_{j=1}^{K}(N-j)\sim NK.
\label{e96}
\end{equation}

Fix a line $\ell(d,t)$ in $\mathcal L_{A,d}^r$ given by the equation
\begin{equation}
y=d+t(x-d).
\label{e97}
\end{equation}
The slope satisfies $t\neq0,1$. 
Suppose that $\ell=\ell(d,t)$ contains $m$ points of $A\times A$, where $r-1\leq m<2r-1$:
\begin{equation}
c_i=d+t(b_i-d),\,b_i\neq d,\,c_i\neq d,
\label{e98}
\end{equation}
\begin{equation}
b_1<b_2<\cdots<b_m.
\label{e99}
\end{equation}

\paragraph{Step 1. Close consecutive pairs.}
Set $K=\lceil8N/r\rceil$. As $t\ne 0$, the sequence $\{c_i\}_{i=1}^m$ is strictly monotone.
Hence
\[
\sum_{i=1}^{m-1}d_I(b_i,b_{i+1})\leq N-1,
\qquad
\sum_{i=1}^{m-1}d_I(c_i,c_{i+1})\leq N-1.
\]
At most $2(N-1)/K\leq r/4$ indices have either gap greater than $K$. Discard also the at most two adjacent pairs containing the pin $(d,d)$. Let
$$
\mathcal J_\ell=\{i\in\{1,\ldots,m-1\}:\,d_I(b_i,b_{i+1})\le K\text{ and }d_I(c_i,c_{i+1})\le K\}.
$$
Since $r\geq16$, the size of $\mathcal J_\ell$ satisfies
\[
|\mathcal J_\ell|\geq m-1-r/4-2\geq r/2.
\]

\paragraph{Step 2. Collinear triples.}
The close pairs define a map
\begin{equation}
(d,t,b_i)\longmapsto(\alpha,\beta_i,\gamma_i),
\qquad i\in\mathcal J_\ell,
\label{e60}
\end{equation}
where
\[
\alpha=(d,d),\qquad \beta_i=(b_i,b_{i+1}),\qquad
\gamma_i=(c_i,c_{i+1}).
\]
The line equation implies
\begin{equation}
\gamma_i-\alpha=t(\beta_i-\alpha).
\label{e61}
\end{equation}
Thus the three points are collinear and $\beta_i,\gamma_i\in P_K$. They are pairwise distinct because $t\neq0,1$ and $\beta_i\neq\alpha$. Their common line is neither horizontal nor vertical, since for any $i$, $b_i\neq d$ and $c_i\neq d$. It is not the diagonal, since $b_i\neq b_{i+1}$. Hence this line belongs to $\mathcal L_{A,d}$.

The map is injective. Its image determines $d$ and $b_i$, and \eqref{e61} determines $t$. The original line and its ordered consecutive points are then fixed. Moreover, the common line of $\beta_i$ and $\gamma_i$ determines its unique diagonal pin. Consequently,
\begin{equation}
\begin{aligned}
r|L_{A,D}^r|
&\lesssim\sum_{\ell\in L_{A,D}^r}|\mathcal J_\ell|\\
&\leq\sum_{\ell'\in L_{A,D}}|\ell'\cap P_K|
\bigl(|\ell'\cap P_K|-1\bigr)\\
&\leq\sum_{\ell'\in L_{A,D}}|\ell'\cap P_K|^2.
\end{aligned}
\label{e62}
\end{equation}
In particular,
\begin{equation}
r|L_A^r|\lesssim\sum_{\ell'\in L_{A,D}}|\ell'\cap P_K|^2.
\label{e63}
\end{equation}

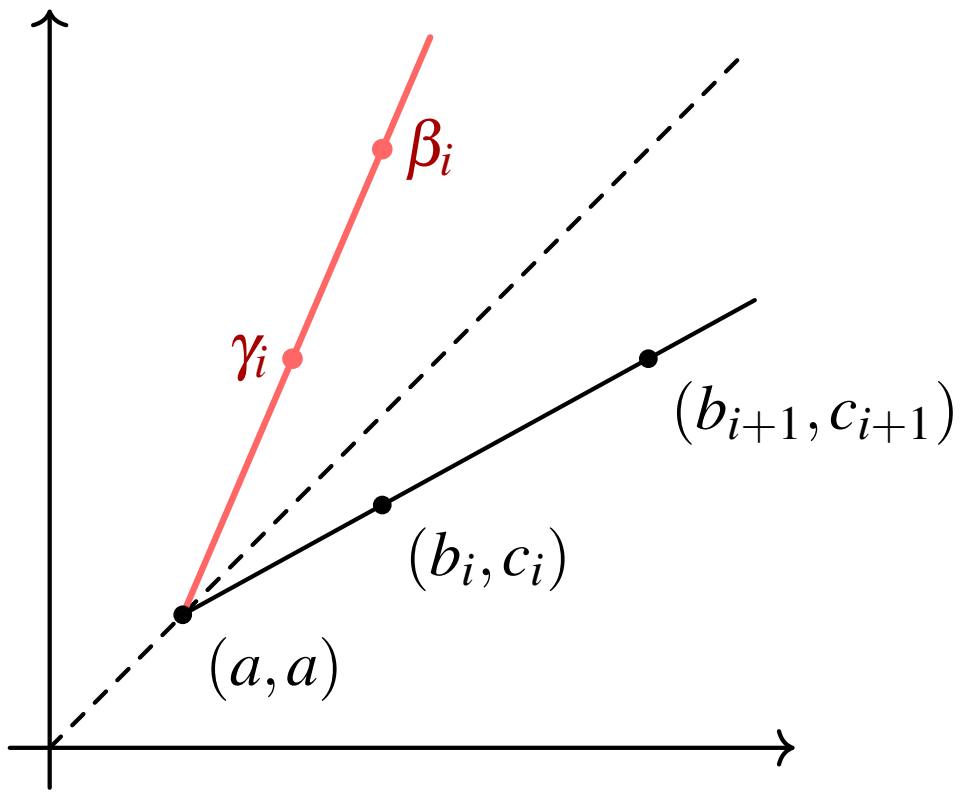
\begin{center}
\begin{tikzpicture}[scale=0.88,line cap=round,line join=round]
    \draw[->,thick] (-0.3,0)--(5.6,0);
    \draw[->,thick] (0,-0.3)--(0,5.55);
    \coordinate (alpha) at (1,1);
    \coordinate (p0) at (2.5,1.825);
    \coordinate (p1) at (4.5,2.925);
    \coordinate (beta) at (2.5,4.5);
    \coordinate (gamma) at (1.825,2.925);
    \draw[dashed,thick] (0,0)--(5.2,5.2);
    \draw[thick] (alpha)--(5.3,3.365);
    \draw[red!60,very thick] (alpha)--(2.86,5.34);
    \fill (alpha) circle (2pt);
    \fill (p0) circle (2pt);
    \fill (p1) circle (2pt);
    \fill[red!60] (beta) circle (2.2pt);
    \fill[red!60] (gamma) circle (2.2pt);
    \node[anchor=north west] at (alpha) {$(a,a)$};
    \node[anchor=north west] at (p0) {$(b_i,c_i)$};
    \node[anchor=north west] at (p1) {$(b_{i+1},c_{i+1})$};
    \node[red!65!black,anchor=west] at (beta) {$\beta_i$};
    \node[red!65!black,anchor=east] at (gamma) {$\gamma_i$};
\end{tikzpicture}
\end{center}

\paragraph{Step 3. A second-moment estimate.}
We prove the bound
$$
\sum_{\ell\in L_{A,D}}|\ell\cap P_K|^2
\lesssim
|D|^{1/2}|P_K|^{3/2}+N|P_K|.
$$

\begin{lm}
\label{t90}
Let $A\subset\mathbb R$ be a finite set. Let $P$ be a subset of
$$
A\times A\setminus\Delta ,
$$
satisfying $|P|\ge|D|$.
Then
\begin{equation}
\sum_{\ell\in L_{A,D}}|\ell\cap P|^2
\lesssim
|D|^{1/2}|P|^{3/2}+N|P|.
\label{e64}
\end{equation}
\end{lm}

\begin{proof}[Proof of the lemma]
The assertion is immediate for $P=\emptyset$. For nonempty $P$ and dyadic $r\geq1$, define
\[
J_{A,D}^r:=\{\ell\in L_{A,D}:r\leq|\ell\cap P|<2r\}.
\]
Since $P$ avoids the diagonal, for fixed $d\in D$, and $p\in P$ there exists only one line $\ell\in L_{A,D}$ satisfying $(d,d),p\in \ell$. Thus
\begin{equation}
\sum_{\ell\in L_A}|\ell\cap P|\leq |D||P|,
\qquad |J_A^r|\leq\frac{|D||P|}{r}.
\label{e65}
\end{equation}
For $r\geq2$, the Szemer\'edi--Trotter theorem also gives
\[
|J_A^r|\lesssim\frac{|P|^2}{r^3}+\frac{|P|}{r}.
\]
We use the bush bound for small $r$ and the Szemer\'edi--Trotter bound for large $r$. Every line has at most $N$ grid points. Set $T=(|P|/|D|)^{1/2}>1$. Summing over dyadic $r$ yields
\[
\begin{aligned}
\sum_{\ell\in L_{A,D}}|\ell\cap P|^2
&\lesssim\sum_{r\leq T}r|D||P|
 +\sum_{T<r\leq N}\left(\frac{|P|^2}{r}+r|P|\right)\\
&\lesssim T|D||P|+\frac{|P|^2}{T}+N|P|\\
&\lesssim |D|^{1/2}|P|^{3/2}+N|P|.
\end{aligned}
\]
The last inequality follows directly from the definition of $T$, including the case $|P|<N$.
\end{proof}

\paragraph{Step 4. Concluding the proof.}
Recall $K\sim N/r$. For the case $|P_K|\ge D$, by \eqref{e96}, \eqref{e63}, and \eqref{e64},
$$
\begin{aligned}
\left|L_{A,D}^r\right|
&\lesssim
\frac{|D|^{1/2}|P_K|^{3/2}}{r}
+
\frac{N|P_K|}{r}
\\
&\lesssim
\frac{N^{3+\frac{\al}{2}}}{r^{5/2}}
+
\frac{N^3}{r^2}.
\end{aligned}
$$
For the case $|P_k|< D$, by \eqref{e96}, we can get that
$$
|D|\gtrsim\frac{N^2}{r}.
$$
Combining with Lemma~\ref{t50}, we can get 
$$
\left|L_{A,D}^r\right|\lesssim\frac{N^4}{r^3}\lesssim\frac{N^{3+\frac{\al}{2}}}{r^{5/2}}.
$$
This concludes the proof.
\end{proof}

\section{Exceptional-set estimates}
\label{s40}

We study the size of exceptional sets related to translates in this section.

\begin{prop}\label{t40}
Let $A\subset\mathbb R$ have cardinality $N\geq2$. Then the exceptional set 
\begin{equation}\label{e22}
S:=\{a\in\bbr:\ |F(a)|<\tf{N(N+1)}2\}.
\end{equation}
satisfies
\[
N^2\lesssim |S|\lesssim N^4.
\]
\end{prop}

\begin{proof}

Suppose $a\in S$, then it is equivalent to say that there exist $(x_1,x_2), (y_1,y_2)\in A\times A$ with $\{x_1,x_2\}\neq\{y_1,y_2\}$ such that 
\begin{equation}\label{e053}
(x_1+x_2-y_1-y_2)a=y_1y_2-x_1x_2.
\end{equation}

We observe first that $(x_1+x_2-y_1-y_2)\neq0$. Otherwise, by \eqref{e053}, $\{x_1,x_2\}=\{y_1,y_2\}$.
As a result,
$a=\frac{y_1y_2-x_1x_2}{x_1+x_2-y_1-y_2}$ has at most $N^4$ possible choices.

For the lower bound, encode unordered pairs by the distinct lines
\[
\mathcal H_A:=\{l_{x,y}:x,y\in A,\ x\leq y\},
\qquad |\mathcal H_A|=M=\frac{N(N+1)}2,
\]
where $l_{x,y}$ is the line given by the map $t\mapsto(x+y)t+xy$. The slope and intercept determine the unordered pair, so these lines are indeed distinct. $a\in S$ if and only if there exists $b\in\bbr$ such that $(a,b)$ is an intersection of two distinct lines in $\mathcal H_A$. Let $\mathcal P$ be the set of intersection points formed by $\mathcal H_A$. Observe that there are at most $N^2$ many points in $\mathcal P$ with a common abscissa. Hence
$$
|\mathcal P|\le N^2 |S|.
$$
It suffices to show that $|\mathcal P|\ges N^4$. Write $d(p)$ for the number of lines through $p\in\mathcal P$. By double counting of the nonparallel lines,
$$
\sum_{p\in\mathcal{P}}\binom{d(p)}2=\sum_{l\in\mathcal H_A}\sum_{l'\in \mathcal H_A}1-\sum_{l\in\mathcal H_A}\sum_{l'\in \mathcal H_A:l'// l}1.
$$
At most $N$ lines have a common slope. So
$$
\sum_{p\in\mathcal P}d(p)^2\ges N^4.
$$
For $r\les N$, the Szemer\'edi--Trotter theorem yields
$$
|\{p\in\mathcal P: d(p)\ge r\}\les \f{M^2}{r^3}+\f{M}{r}\les \f{N^4}{r^3}.
$$
We always have $d(p)\les N$ since for $p=(a,b)\in l_{x,y}$ gives
\[
(x+a)(y+a)=b+a^2
\]
which determines $y$ from $x$ unless its right-hand side is zero; in that case we have $x=-a$ or $y=-a$. As a result, for $R\ge 1$,
$$
\sum_{p\in\mathcal P:d(p)\ge R}d(p)^2\les \f{N^4}{R}.
$$
Taking $R$ a sufficiently large absolute constant guarantees
$$
R^2|\mathcal P|\ge\sum_{p\in\mathcal P:d(p)\le R}d(p)^2\ges N^4.
$$
Consequently, $|\mathcal P|\ges N^4$ as desired.
\end{proof}

We now turn to quantitative bounds at a prescribed threshold. Let $A\subset\bbr$ be finite and nonempty, with $|A|=N$. For fixed $t\in(0,1]$, recall that
\[
S_t:=\{d\in\bbr:|(A+d)(A+d)|\leq N^{2-t}\}.
\]
For sufficiently large $N$ in terms of $t$, this is a subset of the finite set $S$ in \eqref{e22}.

\begin{prop}
\label{t93}
Let $A\subset\bbr$ be finite and nonempty, with $|A|=N$. For every fixed $t\in(0,1]$, there is $N_0(t)$ such that, for $N\geq N_0(t)$,
\[
|S_t|\lesssim_t N^{2-3t}+N^{1-t}\log(2N).
\]
\end{prop}

\begin{proof}
Let $D\subset S_t$ be finite and nonempty. By \eqref{CS}, for each $d\in D$
$$
E_{\times}(A+d)\ge N^{2+t}.
$$
Summing over $d\in D$ yields
$$
\sum_{d\in D}E_{\times}(A+d)\ge |D|N^{2+t}.
$$
We insert this to \eqref{e056} which gives
$$
|D|N^{2+t}
\lesssim N^{8/3}|D|^{2/3}+N^{2}|D|+N^3\log (2N).
$$
As $t>0$, for sufficiently large $N$ in terms of $t$, absorb $N^2|D|$ into the left-hand side. As a result, we have 
$$
|D|N^{2+t}\les N^{8/3}|D|^{2/3}+N^3\log N,
$$
which gives $|D|\les N^{2-3t}+N^{1-t}\log N$. This bound is uniform over all finite subsets $D\subset S_t$, so it also bounds $|S_t|$.
\end{proof}

\section{Examples and counterexamples}
\label{s9}

\subsection{Geometric progressions}\label{exam6.1}

Let $A=\{2^j:1\leq j\leq N\}$. We first verify \eqref{e15}.
Suppose that two distinct unordered pairs from $A$ give the same product
with a nonzero shift $a$, which means that $(a_1+a)(a_2+a)=(a_3+a)(a_4+a)$ where $a_i\in A$ and $\{a_1,\,a_2\}\neq \{a_3,\,a_4\}$. Writing $s_1=a_1+a_2,\,s_2=a_3+a_4$ and $p_1=a_1a_2,\,p_2=a_3a_4$, we have
\begin{equation}\label{e66}
a(s_2-s_1)={p_1-p_2}.
\end{equation}
If $p_1=p_2$, we can get $s_1=s_2$ which will lead to $\{a_1,\,a_2\}=\{a_3,\,a_4\}$. So we can get $p_1\neq p_2$. Combining with the definition of $A$, we can get the following.
\[
|p_1-p_2|\geq \frac12\max\{p_1,p_2\}.
\]
Also, $uv/2\geq u+v-2$ for $u,v\geq2$, so
\[
\frac12\max\{p_1,p_2\}
\geq\max\{s_1,s_2\}-2\geq|s_1-s_2|.
\]
Thus every such collision requires $|a|\geq1$. Consequently,
\[
|(A+a)(A+a)|=\frac{N(N+1)}2\qquad(0<|a|<1).
\]
The same quadratic order holds for every nonzero shift.

\begin{prop}\label{exam1}
Let $A=\{2^j:1\leq j\leq N\}$. Then
\[
|(A+a)(A+a)|\gtrsim N^2\qquad(a\neq0).
\]
\end{prop}
\begin{proof}
If $a$ is irrational, \eqref{e66} rules out all collisions,
and the product set has cardinality $N(N+1)/2$.

Now let $a\in\mathbb Q\setminus\{0\}$ and put $k=v_2(a)$, where
$v_2(2^ku/v)=k$ for odd integers $u,v$. This is the $2$-adic valuation. We use $v_2(zw)=v_2(z)+v_2(w)$ and, when $v_2(z)\neq v_2(w)$,
\[
v_2(z+w)=\min\{v_2(z),v_2(w)\}.
\]
These identities follow by factoring out the smaller power of two. Split the exponents into
\[
J_-=\{1\leq j\leq N:j<k\},\qquad
J_+=\{1\leq j\leq N:j>k\}.
\]
At most one exponent is omitted, so one of these sets has at least $(N-1)/2$ elements. We show that the map
\[
(b,c)\longmapsto P_{b,c}:=(2^b+a)(2^c+a),\qquad b>c,
\]
is injective on pairs from either set.

For $b,c\in J_-$, we have $v_2(P_{b,c})=b+c$. Thus equal products
have the same $b+c$ which means the same $2^{b+c}$, and expanding them gives the same value of
$2^b+2^c$. Uniqueness of binary expansion then determines $(b,c)$.
For $b,c\in J_+$,
\[
v_2(P_{b,c}-a^2)
=v_2\bigl(2^{b+c}+a2^c(2^{b-c}+1)\bigr)=k+c.
\]
Indeed, the two summands have different valuations $b+c$ and $k+c$.
Hence $P_{b,c}$ determines $c$, and then $b$, since $2^c+a\neq0$.
Choosing pairs from the larger of $J_-$ and $J_+$ therefore gives
$\gtrsim N^2$ distinct products for sufficiently large $N$; bounded $N$
is immediate. This also covers shifts $a=-2^k$ without a separate argument.
\end{proof}

\subsection{Algebraic integers}
In this subsection we provide a proof of Theorem~\ref{t30} using the ingredients of \cite{B26} and \cite{R26}. We start with a number field $K$ with degree $n=[K:\mathbb Q]<\infty$. Then there are exactly $n$ distinct embeddings $\sigma_i:K\rightarrow\bbc$, $i=1,2,\cdots,n$. 
Note that each $\sigma_i$ is an injection that preserves rational numbers. 
The number field $K$ is called totally real if $\sigma_i(K)\subset\bbr$ for all $i=1,2,\cdots,n$. In this section we always consider totally real number fields. 

A number $x\in K$ is an algebraic integer if it is a root of some monic polynomial with integer coefficients. Denote the family of algebraic integers in $K$ by $\mathcal O_K$. Clearly, $\bbz\subset\mathcal O_K$. The monic condition forces that $\mathcal O_K\cap \mathbb Q=\bbz$. The following classical lemma will be useful for us.
\begin{lm}
    $\mathcal O_K$ is a ring containing $\bbz$.
\end{lm}
Then we denote the unit group of $\mathcal O_K$ by $\mathcal O_K^{\times}$. Note that $\sigma_i(u)\ne 0$ for any $u\in\mathcal O_K^{\times}$. Moreover, $u\in\mathcal O_K^{\times}$ if and only if $Nm_{K/\mathbb Q}(u)=\prod_{i=1}^n\sigma_i(u)=\pm1$. We introduce the following two maps.
$$
\iota:\mathcal O_K\rightarrow\bbrn,\,\iota(x):=(\sigma_1(x),\sigma_2(x),\cdots,\sigma_n(x)).
$$
$$
\ell:\mathcal O_K^{\times}\rightarrow\bbr^n,\,\ell(u):=(\log|\sigma_1(u)|,\log|\sigma_2(u)|,\cdots,\log|\sigma_n(u)|).
$$
Let $\La=\iota(\mathcal O_K)$ and $\La_{\log}=\ell(\mathcal O_K^{\times})$.
With the help of $\bbz$-basis for $\mathcal O_K$ and $\mathcal O_K^{\times}$, we can see that $\La$ is a full-rank lattice in $\bbrn$ and $\La_{\log}$ is an $(n-1)$-dimensional lattice in the hyperplane $H=\{y\in\bbrn:y_1+\cdots+y_n=0\}$. Under the above two maps, we can define the following additive and multiplicative boxes.
\begin{defn}
    For $X,Y\ge 1$, define
    $$
    B_K^{+}(X):=\{x\in\mathcal O_K:\|\iota(x)\|_{\infty}\le X\},
    $$
    $$
    B_K^{\times}(Y):=\{u\in\mathcal O_K^{\times}:\|\ell(u)\|_{\infty}\le Y\}.
    $$
\end{defn}
As each $\sigma_i$ is injective, it holds that $|B_K^+(X)|=|\La\cap[-X,X]^n|$ and $|B_K^{\times}(Y)|=2|\La_{\log}\cap[-Y,Y]^n|$. The following lemma describes how these boxes change under addition and multiplication.
\begin{lm}\label{t92}
    For $X_1,X_2,Y_1,Y_2\ge1$, 
    $$
    B_K^+(X_1)+B_K^+(X_2)\subset B_K^+(X_1+X_2),
    $$
    $$
    B_K^+(X_1)B_K^+(X_2)\subset B_K^+(X_1X_2),
    $$
    $$
    B_K^{\times}(Y_1)B_K^{\times}(Y_2)\subset B_K^{\times}(Y_1+Y_2),\quad B_K^{\times}(Y_1)^{-1}=B_K^{\times}(Y_1),
    $$
    $$
    B_K^{\times}(Y_1)\subset B_K^+(e^{Y_1}).
    $$
\end{lm}

The following theorem is guaranteed by \cite{B26} and \cite{M76}.
\begin{thm}\label{BC}
    There exists an absolute constant $B\ge 2$ and a sequence of totally real number fields $K$ with $n:=[K:\mathbb Q]\rightarrow\infty$ such that
    $$
    \Delta_K\le B^n,
    $$
    where $\Delta_K$ is the discriminant of the number field $K$. 
    For any totally real number field $K$ of degree $n\ge2$ and any $X,Y\ge1$, we have
    $$
    \f{X^n}{\sqrt{\Delta_K}}\le |B_K^+(X)| \le (2X+1)^n,
    $$
    $$
    \f{Y^{n-1}}{\sqrt{n}\Delta_K}\le |B_K^{\times}(Y)| \le 10(5Y+1)^{n-1}.
    $$
\end{thm}

We also need the following lemma proven in \cite{B26}.
\begin{lm}\label{N1L}
    Fix $\ep_0=1/8$. Let $T>0$ and let $P_T\subset\mathcal O_K$ satisfy
    $$
    (1-\ep_0)T\le\sigma_i(p)\le(1+\ep_0)T,
    $$
    for any $i=1,2,\dots,n$ and $p\in P_T$. Then for any $G\subset\mathcal O_K^{\times}$, the multiplicative map
    $$
    G\times P_T\rightarrow GP_T:(g,p)\mapsto gp
    $$
    is injective. In particular, for finite $G$ and $P_T$ as above, $|GP_T|=|G||P_T|$.
\end{lm}

With the above tools in hand, we can start to prove Theorem~\ref{t30}.
\begin{proof}[Proof of Theorem~\ref{t30}]
Fix a totally real number field $K$ of degree $n\ge2$ from the sequence in Theorem~\ref{BC}. Let $Y\ge1$, $X,Z\in\bbz_+$ be parameters to be chosen later satisfying $\ep_0X\ge1$, $\ep_0 Z\ge 1$ and $Xe^{-2Y}/8\le Z\le Xe^{-2Y}/4$. Define
$$
G:=B_K^{\times}(Y),
$$
$$
P:=X+B_K^+(\ep_0 X),\,Q:=Z+B_K^+(\ep_0Z),
$$
$$
A:=GP,\, D:=GQ.
$$
As $X,Z\in\bbz$, we have $\sigma_i(X)=X$ and $\sigma_i(Z)=Z$ for any $i$. Therefore, for any $p\in P$ and $q\in Q$, it holds that $\sigma_i(p)\in [(1-\ep_0)X,(1+\ep_0)X]$ and $\sigma_i(q)\in [(1-\ep_0)Z,(1+\ep_0)Z]$. So Lemma~\ref{N1L}
implies 
$$
N:=|A|=|G||P|,\,|D|=|G||Q|.
$$
Notice that $|P|=|B_K^+(\ep_0X)|$ and $|Q|=|B_K^+(\ep_0Z)|$. As a result of Theorem~\ref{BC}, there is a large absolute constant $C\gg1$ depending on $B$ which may change from line to line such that
\begin{equation}\label{e107}
    C^{-n}{Y^{n-1}X^{n}}\le N\le C^n  Y^{n-1}X^{n},
\end{equation}
\begin{equation}\label{e108}
    |D|=N\f{|Q|}{|P|}\ge N\f{(\ep_0 Z)^n}{\sqrt{\Delta_K}(2\ep_0 X+1)^n}\ge C^{-n}Ne^{-2nY},
\end{equation}
where we have used the constraint $Z\ge Xe^{-2Y}/8$.

To estimate $M_D(A)$, we estimate $|(A+D)(A+D)|$. Notice that for $a_j=g_jp_j\in A$ and $d_j=u_j q_j\in D$, where $g_j,\,u_j\in G$, $p_j\in P$ and $q_j\in Q$, $j=1,2$,
$$
(a_1+d_1)(a_2+d_2)=g_1g_2(p_1p_2+\f{u_1}{g_1}p_2q_1+\f{u_2}{g_2}p_1q_2+\f{u_1u_2}{g_1g_2}q_1 q_2).
$$
Combing with $B_K^{\times}(Y_1)^{-1}=B_K^{\times}(Y_1)$, we can get that $(A+D)(A+D)\subset G^2(P^2+2G^2 PQ+G^4 Q^2)$. 
Then, by Lemma~\ref{t92} and our condition on $Z$, we conclude that 
$$
(A+D)(A+D)\subset B_K^{\times}(2Y)B_K^+(5X^2).
$$
An again application of Theorem~\ref{BC} yields
\begin{equation}\label{e109}
    |(A+D)(A+D)|\le C^n X^{2n}{Y^{n-1}}\le C^{2n} Y^{-(n-1)}N^2.
\end{equation}

Let  $c_*\in(0,1/4)$ be a sufficiently small to be determined constant independent of $\al$. 
Now for a fixed $\al\in[0,1)$, by \eqref{e108} and \eqref{e109}, to prove Theorem~\ref{t30} we only need to choose $X,Y$ such that
$$
C^{-n}e^{-2nY}\ge (C^{-n}Y^{n-1}X^n)^{-(1-\alpha)}\ge N^{-(1-\al)},
$$
and 
$$
C^{2n} Y^{-(n-1)}\le (C^{n}Y^{n-1}X^n)^{-c_*(1-\alpha)}\le N^{-c_*(1-\al)},
$$
which follows from
\begin{equation}\label{e210}
    Y+(\frac{1-\al}{2}+1)\log C\le \f{1-\al}{2}\log X
\end{equation}
and 
\begin{equation}\label{e211}
    c_*(1-\al)\log X+[\frac
{c_*(1-\al)}{2}+2]\log C\le \frac 14\log Y
\end{equation}
respectively.

To achieve \eqref{e210} and \eqref{e211}, we take $Y\approx e^{100C}$, by which we mean $|Y-e^{100C}|\le 5$, and
$X$ such that $\f{1-\al}{2}\log X\approx 2e^{100C}$. Then \eqref{e210} directly holds and \eqref{e211} reduces to
$$
5c_*e^{100C}\le 25C.
$$
This can be guaranteed by taking $c_*=e^{-200C}$. Our choices of $X$ and $Y$ also ensure that $Xe^{-2Y}\ge e^{2C}$ so we can choose $Z$ satisfying the desired conditions. Moreover, $N^2\ge C^{-n} X^{2n}Y^{2n-2}$ will tend to infinity as $n\rightarrow\infty$.

We complete the proof of Theorem~\ref{t30}.
\end{proof}

\end{document}